\input amstex
\input psfig.sty

\documentstyle{amsppt}

\NoRunningHeads
\magnification\magstep 1
\baselineskip 15pt
\pagewidth{6.4 truein}
\pageheight{8.6 truein}

\TagsOnRight

\catcode`\@=11
\redefine\logo@{}
\catcode`\@=13

\def\qed{\hfill\ \vbox{\hrule\hbox{\vrule\kern4pt\vbox{\kern4pt{}
\kern4pt}\kern4pt\vrule}\hrule}}

\hfill 8/26/98%

\topmatter
\title 
Control of Nonlinear Underactuated Systems
\endtitle

\author  {\ }\endauthor

\affil
{\smc David Auckly}\\
Department of Mathematics, \\
Kansas State University, Manhattan, KS 66506-2602, USA\\
email: dav\@math.ksu.edu\\
\\
{\smc Lev Kapitanski} \\
Department of Mathematics,   \\
Kansas State University, Manhattan, KS 66506-2602, USA\\
email: levkapit\@math.ksu.edu \\
\\
{\smc Warren White} \\
Department of Mechanical and Nuclear Engineering\\
Kansas State University, Manhattan, KS 66506-5106, USA\\
email: wnw\@ksu.edu \\
\\
\endaffil


\endtopmatter


In this paper we introduce a new method to design control laws for
non-linear underactuated systems. Our method will often 
work in situations
where standard control techniques fail. Even when standard
techniques apply, we believe that our approach will achieve better
performance. This is because we produce an infinite dimensional
family of control laws, whereas most control techniques only
produce a finite dimensional family of control laws. We will describe
the problem and our solution invariantly, using
differential geometry and in local coordinates. We hope that
this paper will be useful for both mathematicians and engineers.
We include an abstract example of a system which is open-loop
unstable and cannot be stabilized using any linear control law, and
demonstrate that our method produces a stabilizing control law.
We also apply our method to the inverted pendulum cart to compare it with
some standard control techniques.

An important problem in control theory is how to modify a system of ordinary
differential equations so that solutions to the equations satisfy
some property. In order to describe the class of differential equations
that we consider, let $\,Q\,$ denote the configuration space.
The configuration space is a finite dimensional manifold that
represents, for example, every possible position of a mechanical
system. Let $g\in\Gamma(T^\ast Q\otimes T^\ast Q)$ be a metric.
This is the mass matrix in a mechanical system.
Let $c,f:TQ\to TQ$ be fiber-preserving maps. These maps are not
assumed to be linear, but we do assume that $c$ is odd, i.e.,
$c(-X)=-c(X)$. In a mechanical system $c$ will represent the inherent
dissipation and $f$ will represent the input forces.
Finally, let $V:Q\to{\Bbb R}$. This will represent the potential
energy of a mechanical system. The
differential equation that we consider in this paper is
$$ \nabla_{\dot\gamma}\dot\gamma 
+ c(\dot\gamma) +\ grad_\gamma V
   =f(\dot\gamma).  \tag1 
$$
Let $P\in\Gamma(T^\ast Q\otimes TQ)$ be a $g$-orthogonal projection. 
We consider the situation where a constaint $P(f) = 0$ is imposed.
A system is called underactuated if $P \neq 0$. 
A mechanical system is underactuated if we require
that control forces are zero in certain directions.
Many problems in chemical, electrical and mechanical engineering
may be formulated as follows: find a function $f$ with $P(f)=0$ so
that solutions to Equation (1) have some specified property.
For example, solutions with initial conditions in some region
will remain close to some path, or as a different example,
some point in $TQ$ will be an asymptotically stable equilibrium.

In local coordinates $\,\underline{x}=(x^1,\dots,x^n)\,$ 
Equation (1) reads
$$ \ddot x^k + \Gamma^k_{ij}(\underline x) \dot x^i\dot x^j
   + c^k(\underline{x}, \dot{\underline x})
   + g^{ik} \frac{\partial V}{\partial x^i}
   = f^k(\underline{x},\dot{\underline x}), \tag 2
$$
where $\,\Gamma^k_{ij}(\underline{x})\,$ are the
Christoffel symbols (of the Levi-Civita connection) 
associated to the metric $\,g$, \cite{Hicks, 1965}, 
and $g^{ik}$ is the inverse matrix to $g_{kj}$. 
In  Equation (2) and throughout the rest of this paper
we are using the
summation convention, that repeated indices are summed from
1 to $n$. 
As we have mentioned,  we assume that
$c^k(\underline{x},-\dot{\underline x})
  =-c^k(\underline x,\dot{\underline x})$
and the projection of $f^k$ to
some specified set of directions must be zero.

The main question addressed in this paper is how to find a function, 
$f$, with
$P(f)=0$, so that solutions to Equation (1) satisfy some preassigned
conditions.
Our approach to this question is to find functions 
$\widehat g$, $\widehat c$, 
$\widehat V$ and $f$ so that
solutions to Equation (1) are automatically solutions to
$$ \widehat \nabla_{\dot\gamma} \dot\gamma 
+ \widehat c(\dot\gamma) +\ \widehat  {grad}_\gamma
   \widehat V=0. $$
This will clearly be the case if
$$ f(X)\equiv \nabla_X X-\widehat \nabla_X X+\ grad_\gamma V
   -\ \widehat {grad}_\gamma \widehat V +c(X)-\widehat c(X),
   \tag3 
$$
for every vector field $\,X$.
The condition $P(f)=0$ then becomes a system of nonlinear partial
differential equations for $\widehat g$, $\widehat c$, and 
$\widehat V$.
Notice that constant multiples of $g$, $c$, and $V$ satisfy
$\,P(f)=0$ even when $\,P$ has full rank. Thus, one would expect many
solutions when $\,P$ does not have full rank. Separating
$P(f)=0$ into terms which are quadratic in the velocity, independent
of the velocity or
odd functions of the velocity gives
$$ 
P(\nabla_X X-\widehat \nabla_X X)=0, \tag4.1 
$$
$$ 
P(grad_\gamma V-\ \widehat  {grad}_\gamma \widehat V)=0, \tag4.2 
$$
$$ 
P(c(X)-\widehat c(X))=0. \tag4.3 
$$
We will look for solutions to these matching equations with
$\widehat g$ non-degenerate so that $g(X,Y)=\widehat g(\lambda X,Y)$ with
$\lambda\in\Gamma(T^\ast Q\otimes TQ)$. It is clear that $\lambda$
has to be
$g$ self-adjoint, i.e., $g(\lambda X,Y)=g(X,\lambda Y)$. 
We will derive a linear system of partial
differential equations for $\lambda$ which must be
satisfied if $\,\widehat g$ is to solve Equation (4.1).

To derive the equations for $\lambda$ we will use the relation
between the connection, Lie bracket (commutator) and metric. 
In fact, we only derive equations for $\lambda\big|_{\text{Im}\ P}$. 

It is known, \cite{Hicks}, that the covariant derivative 
$\,\nabla\,$ compatible with the metric $\,g\,$ 
is determined uniquely as a bilinear  
operator which associates to any pair of vectors
$X$ and $Y$ a third vector $\nabla_XY$ so that the following 
equations are satisfied:
$$ 
X(g(Y,Z))=g(\nabla_XY,Z)+g(Y,\nabla_XZ),
   \tag5 
$$
$$
\nabla_XY-\nabla_YX=[X,Y],
   \tag6 
$$
where $\,[X,Y]\,$ is the Lie bracket (commutator) of $X$ and $Y$. 
Using the above properties of the covariant derivative we get:
$$ \align
   2g(\nabla_XY,Z) =
    & Xg(Y,Z) +Yg(Z,X)-Zg(X,Y) \\
    &+g([X,Y],Z)+g([Z,X],Y)-g([Y,Z],X). \tag7 \endalign
    $$
This may be solved for $\nabla_XY$ since $g$ is non-degenerate.

\proclaim{Proposition 1}
If $g=\widehat g\lambda$ and $\widehat g$ 
satisfies the matching condition
$P(\nabla_X X-\widehat \nabla_X X)=0$, then $\lambda$ satisfies:
$$ 
\nabla g\lambda\big|_{\text{Im}\ P^{\otimes 2}}=0,
   \tag8 
$$
or, equivalently,
$$ 
g(\nabla_Z\lambda PX, PX)-g(\lambda PX, \nabla_Z PX)=0.
   \tag9 
$$
\endproclaim

\demo{Proof}
We begin by polarizing the matching equation
$$ 
\align
   P(\nabla_XY-\widehat \nabla_XY)
    & =\frac12 P(\nabla_XY+\nabla_YX
                      -\widehat \nabla_XY-\widehat \nabla_YX) \\
    & = \frac12 \left[ P(\nabla_{X+Y}(X+Y)
        -\widehat \nabla_{X+Y}(X+Y)) \right.\\
    &\qquad \left. -P(\nabla_XX-\widehat \nabla_XX)
         -P(\nabla_YY-\widehat \nabla_YY) \right]\\
    &=0.
    \endalign 
$$
The first line is true because the covariant derivatives are
torsion free (Equation (6)). The second line is true because the
covariant derivatives are bilinear (Equation (7)). Now,
$$ \align
 0 &= 2g(P(\nabla_{\lambda PX}Z-\widehat \nabla_{\lambda PX} Z),X) \\
   &= 2g(\nabla_{\lambda PX} Z-\widehat \nabla_{\lambda PX}Z,PX  ) \\
   &= 2g(\nabla_{\lambda PX} Z,PX)
      -2\widehat g(\widehat \nabla_{\lambda PX}Z,\lambda PX) \\
   &=\lambda PX g(Z,PX)+Zg(PX,\lambda PX) -PXg(\lambda PX,Z) \\ 
   &\qquad +  g([\lambda PX,Z],PX) + g([PX,\lambda PX],Z)
       -g([Z,PX],\lambda PX) \\
   &\qquad -\lambda PX \widehat g(Z,\lambda PX)
      -Z\widehat g(\lambda PX, \lambda PX)
       +\lambda PX \widehat g(\lambda PX,Z) \\
   &\qquad- \widehat g([\lambda PX, Z],\lambda PX)
         -\widehat g([\lambda PX,\lambda PX],Z)
       +\widehat g([Z,\lambda PX],\lambda PX) \\
   &= \lambda PX g(Z,PX)-PX g(\lambda PX,Z)+g([PX,\lambda PX],Z)\\
   &\qquad -g([Z,PX],\lambda PX) +g([Z,\lambda PX],PX)\\
   &= g(\nabla_{\lambda PX}Z,PX) +g(Z,\nabla_{\lambda PX}PX) \\
   & \qquad -g(\lambda PX, \nabla_{PX}Z)-g(Z,\nabla_{PX}\lambda PX)\\
   & \qquad +g([Z,\lambda PX],PX) -g([PX,\lambda PX],Z)
     -g([Z,PX],\lambda PX) \\
   &= g(\nabla_Z(\lambda PX),PX)-g(\nabla_ZPX,\lambda PX) \\
   &= g((\nabla_Z\lambda)(PX),(PX)) \\
   &= (\nabla_Z\, g\lambda)(PX,PX).
   \endalign $$
The covariant derivative is extended to all tensors by requiring every
reasonable product rule to hold. For example,
$\, 
Z(\lambda(X)) =\nabla_Z(\lambda(X))
=(\nabla_Z\lambda)(X)+\lambda(\nabla_ZX)$
and
$\,
X(g(Y,Z)) =\nabla_X(g(Y,Z)) 
= (\nabla_Xg)(Y,Z) +g(\nabla_XY,Z) +g(Y,\nabla_X Z)$.
We used these two product rules in the last two lines and Equations (5),
(6) and (7) in the previous lines. \qed
\enddemo

Now that we have equations for $\lambda\big|_{\text{Im}\ P}$, 
we will derive equations
for $\widehat g$ in terms of $\lambda\big|_{\text{Im}\ P}$.

\proclaim{Proposition 2}
If $g=\widehat g\lambda$, then
$$ \lambda PX\, \widehat g(Z,Z) +2\widehat g([Z,\lambda PX],Z)
   =2Zg(PX,Z)-2g(PX,\nabla_Z Z)
   \tag10 $$
\endproclaim

\demo{Proof}
$$ 
\align
  \lambda PX\, &\widehat g(Z,Z)+2\,\widehat g([Z,\lambda PX],Z) \\
  &=2\,\widehat g(\widehat \nabla_{\lambda  PX} Z,Z)
     +2\,\widehat g([Z,\lambda PX],Z)\\
  &= 2\,\widehat g(\widehat \nabla_Z \lambda PX,Z)\\
  &= 2\,Zg(PX,Z)-2g(PX,\nabla_ZZ). \tag"\qed"
  \endalign 
$$
\enddemo

The matching equation for the potential energy may be
expressed as a linear partial differential equation using $\lambda$.

\proclaim{Proposition 3}
If the matching Equations (4) are satisfied and $g=\widehat g\lambda$,
then
$$
\lambda PX(\widehat V)=PX(V).
   \tag11 
$$
\endproclaim

\demo{Proof} 
$\quad 
\lambda PX(\widehat V) \,=\, (d\widehat V)(\lambda PX)
 \,=\,\widehat g(\widehat  {grad}_\gamma \widehat V,\lambda PX)\,
   = \,g(\widehat  {grad}_\gamma \widehat V, PX)\,$ 
\newline\noindent
 $=\, g(P\ \widehat  {grad}_\gamma \widehat V,X)
   \,=\, g(P\ grad_\gamma V,X)=PX(V)$. \qed
\enddemo

By solving the matching equations, we find every control law that will
result in dynamical behavior equivalent to a system of the form: 
$$ \widehat \nabla_{\dot\gamma} \dot\gamma + \widehat c(\dot\gamma) +\ \widehat {grad}_\gamma
   \widehat V=0. $$
This is a very large collection of vector fields on $TQ$. Thus far, we
have shown that every solution to the matching equations may be found
by first solving one system of linear partial differential equations
and then solving a different system of linear partial differential
equations.

The previous three propositions suggest a natural approach for solving the
matching equations. First, solve Equation (9) for $\lambda PX$,
then solve Equation (10) for $\widehat g$, and Equation (11) 
for $\widehat V$.
Finally solve the algebraic Equation (4.3) for $\widehat c$. 
Every solution  to
the matching equations can be found in this way. What is not
clear is whether every set of functions generated in this way 
is a solution
to the matching equations. The problem is that once we
have a solution to the $\lambda$-equation, (Equation (9)), 
and then
a solution to Equation (10), it is not clear
that $\lambda$ can be extended from $\text{Im}\ P$ so that the 
condition 
$g=\widehat g\lambda$ holds. We will address this question for 
general $P$ in a future paper.
We will next show that $g=\widehat g\lambda$ holds if it is true at
the initial conditions, provided that the rank of $P$ is one.

\proclaim{Proposition 4}
Let $X$ be a non-zero vector field in an open set $U$,
generating $\text{Im}\ P$, and let $\lambda PX$ be a non-vanishing
solution of the $\lambda$ Equation (9) in $U$. Let $\Sigma  $ be
a codimension one non-characteristic hypersurface in $U$.
Assume that $\widehat g$ is defined, non-degenerate, and
satisfies $\widehat g\lambda PX=gPX$ on $\Sigma$, 
and $\widehat V$ is defined
on $\Sigma $. Then there is a unique solution to the matching
Equations (4.1) and (4.2), and a unique extension of $\lambda$
away from $\text{Im}\ P$ in a neighborhood of $\Sigma  $.
Furthermore, $\widehat g\lambda=g$ everywhere in the neighborhood.
\endproclaim

\demo{Proof} Notice that $\Sigma$ is non-characteristic
for the $\widehat g$-Equation (10) if it is non-characteristic for the
$\widehat V$-Equation (11), and {\it visa versa\/}. 
The functions $\widehat g$ and $\widehat V$
may be found in a neighborhood of $\Sigma$ 
using the method of characteristics.
The main content of this proposition
is that $g PX=\widehat g\lambda PX$ everywhere in the neighborhood.
After polarization, the $\widehat g$-Equation (10) reads:
$$ 
\align
 & \lambda PX\, \widehat g(\lambda PX,Z)
+\widehat g(\lambda PX,[Z,\lambda PX])\\
 &\qquad =Zg(\lambda PX,PX)+\lambda PX g(PX,Z)
    -g(\nabla_Z\lambda PX,PX)-g(\nabla_{\lambda PX}Z,PX) \\
 &\qquad = \lambda PX\, g(PX,Z)+g(PX,[Z,\lambda PX])
    +g(\lambda PX,\nabla_ZPX)-g(\nabla_Z\lambda PX,PX) .
\endalign 
$$
This last line follows after applying the product rule to the first
term of the previous line and the torsion-free condition (6)
to the last term of the previous line. In view of the 
$\lambda$-Equation (9), 
it is clear that $\widehat g\lambda PX=gPX$ is a solution
of this equation. This solution is unique in a neighborhood
of $\Sigma  $. Finally, since $\widehat g$ is non-degenerate on
$\Sigma  $, it is non-degenerate in a neighborhood of
$\Sigma  $, so we may extend $\lambda$ to the full tangent
space by $\lambda=(\widehat g)^{-1}g$. \qed
\enddemo
\medskip
\centerline{ EXAMPLE}

In order to demonstrate and clarify our method, we will apply it to an
inverted pendulum cart as an example.
\medskip

\hskip 100bp\psfig{file=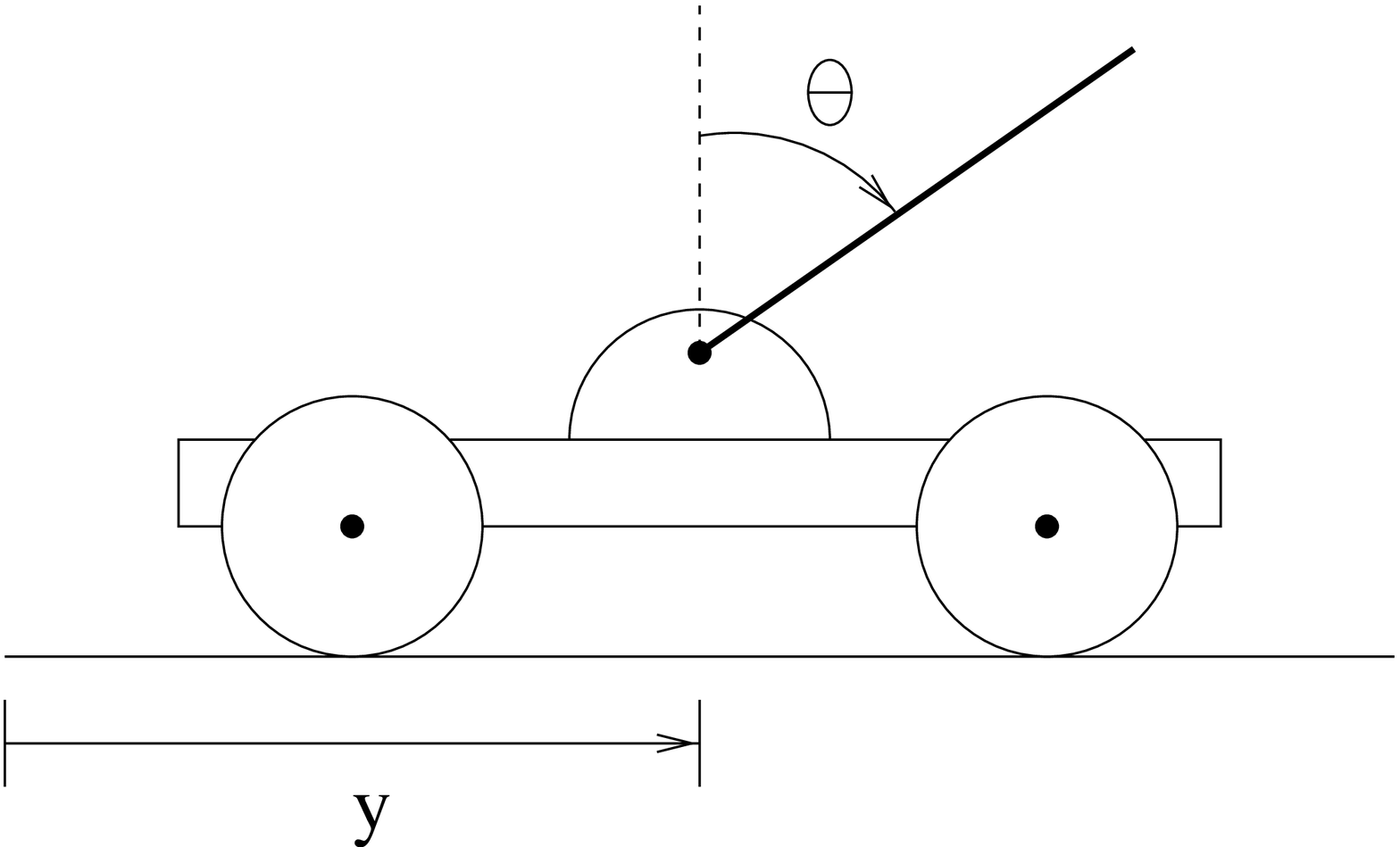,width=3.1truein,height=2truein}%

\centerline{Figure 1}

\vfill
\pagebreak

The (mass) metric for the cart system  depicted in Figure~1 is:
$$ 
g  
=(M+m) dy^2+ 2 m\ell \cos(\theta)d\theta dy
   +(m\ell^2+I)d\theta^2,
$$
and the potential energy is:
$$ 
V=m\text{\rm g}\ell\cos(\theta). 
$$
Here
\roster
\item"" $M\;\,$ is the mass of the base of the cart,

\item"" $m\;\;$ is the mass of the pendulum,

\item"" $\ell\;\quad\!$ is the length from hinge to the center
of mass of the pendulum,

\item"" $I\;\quad\!$ is the moment of inertial about the
center of mass, and 

\item"" $\text{\rm g}\;\quad\!$ is the acceleration due to gravity.
\endroster

By a change of length, time and mass scales, the metric and 
the potential energy 
may be transformed into
$$  
g = dx^2
+ 2 b\cos(\theta)d\theta dx+d\theta^2,
      \quad V=\cos(\theta),
   \tag12 
$$
where $b=m\ell(M+m)^{-\frac12}(m\ell^2+I)^{-\frac12}\in (0,1)$
is a dimensionless parameter.

The problem is to find the force applied to the base of the cart
as a function of the state of the cart so that the origin
$(\theta=0, x=0)$ will be an asymptotically stable equilibrium of
the resulting dynamical system.

In this application, we may only directly apply force in the
$x$ direction, so $P$ is projection
onto the $\frac{\partial}{\partial\theta}$ direction, i.e.,
$$ 
P=(b\cos(\theta) \,dx+d\theta)\otimes\frac{\partial}{\partial\theta}.
$$
Projecting system (1) onto the coordinate directions gives:
$$
\aligned 
P\left(\nabla_{\dot\gamma}\dot\gamma 
+ c(\dot\gamma) +\ grad_\gamma V
   \right) & =0, \\
g\left( \frac{\partial}{\partial x},
   \nabla_{\dot\gamma}\dot\gamma + c(\dot\gamma) +\ grad_\gamma V\right)
   &=u\equiv g\left( \frac{\partial}{\partial x},f \right).
\endaligned \tag 13
$$
We will assume that there is no
inherent dissapation: $\,c\equiv 0$. Then,
in local coordinates, the equations read:
$$
\aligned 
\ddot\theta + b\cos(\theta)\ddot x -\sin(\theta) &= 0 \\    
b\cos(\theta)\ddot\theta + \ddot x -b\sin(\theta)\dot\theta^2 &= u.
\endaligned
\tag 14
$$
Equation (3) expresses our solution to this problem in terms of
functions, $\widehat g_{ij}$, $\widehat V$ and $\widehat c_i$, 
$i,j=1$ or $2$,
which satisfy the Equation (4).
In local coordinates, these matching equations read: 
$$ 
   (\widehat g_{22}-b \cos(\theta)\widehat g_{12}) 
  \frac{\partial \widehat g_{11}}{\partial\theta}
  +
  (b\cos(\theta)\widehat g_{11}-\widehat g_{12})
  \left[ 2\frac{\partial \widehat g_{12}}{\partial\theta}
    -\frac{\partial \widehat g_{11}}{\partial x} \right] =0,
  $$
$$ 
  (\widehat g_{22}-b\cos(\theta)\widehat g_{12})
  \left[ 2\frac{\partial \widehat g_{12}}{\partial x}
    -\frac{\partial \widehat g_{22}}{\partial\theta} \right] 
   +
  (b\cos(\theta)\widehat g_{11}-\widehat g_{12})
  \frac{\partial \widehat g_{22}}{\partial x} =0,
  $$
$$ 
  (\widehat g_{22}-b\cos(\theta)\widehat g_{12})
  \frac{\partial \widehat g_{11}}{\partial x}
  +
  (b\cos(\theta)\widehat g_{11}-\widehat g_{12})
  \frac{\partial \widehat g_{22}}{\partial\theta} =0
  $$
$$\qquad
     (\widehat g_{22}-b\cos(\theta)\widehat g_{12})
     \frac{\partial \widehat V}{\partial\theta} 
+ (b\cos(\theta)\widehat g_{11}-\widehat g_{12})
     \frac{\partial \widehat V}{\partial x}     
=-\sin(\theta)(\widehat g_{11}\widehat g_{22}-(\widehat g_{12})^2),
  $$
and
$$ 
\widehat c_1+b\cos(\theta)\widehat c_2=0 .    
$$
We have slightly simplified these equations by multiplying 
by the determinant
of $\widehat g$.  
(The terms on the right hand side 
of each equation would be non-zero if the cart was on a hill, and
the hinge was rusty. 
Therefore, multiplying by the determinant of $\,\widehat g\,$ 
would not be a great simplification.) 
These equations are messy, but the same functions may be found
by solving the equations  in Propositions 1, 2, and 3 without ever
writing out the matching equations. Take 
$\,X=\frac{\partial}{\partial\theta}\,$
in Proposition 1 (so $PX=\frac{\partial}{\partial\theta}$) 
and let
$$ 
\lambda\left( \frac{\partial}{\partial\theta}\right)
  =\sigma\frac{\partial}{\partial\theta}
  +\mu\frac{\partial}{\partial x}.
  $$
Using the bilinearity of the covariant derivative and the product rule,
the $\lambda$-Equation (9) becomes:
$$  
g\left( \frac{\partial}{\partial\theta},
   \frac{\partial}{\partial\theta}  \right)
   \frac{\partial\sigma}{\partial\theta} 
   +g\left( \frac{\partial}{\partial x},
   \frac{\partial}{\partial\theta}  \right)
   \frac{\partial\mu}{\partial\theta} 
    + g\left( \nabla_{\frac{\partial}{\partial\theta}}
 \frac{\partial}{\partial x}, \frac{\partial}{\partial\theta} \right)\mu
    - g\left( \frac{\partial}{\partial x}, 
 \nabla_{\frac{\partial}{\partial\theta}} \frac{\partial}{\partial\theta}
    \right) \mu=0
  $$
and
$$ g\left( \frac{\partial}{\partial\theta},
   \frac{\partial}{\partial\theta}  \right)
   \frac{\partial\sigma}{\partial x}
 +g\left( \frac{\partial}{\partial x},
   \frac{\partial}{\partial\theta}  \right)
   \frac{\partial\mu}{\partial x} 
   + g\left( \nabla_{\frac{\partial}{\partial x}}
 \frac{\partial}{\partial x}, \frac{\partial}{\partial\theta} \right)\mu
    -g\left( \frac{\partial}{\partial x},
 \nabla_{\frac{\partial}{\partial x}} \frac{\partial}{\partial\theta}
    \right) \mu=0
  $$
Using Equation (7), it follows that:
$$ \align
 \frac{\partial\sigma}{\partial\theta}
   +b\cos(\theta) \frac{\partial\mu}{\partial\theta}
   +b\sin(\theta)\mu &=0, \\
 \frac{\partial\sigma}{\partial x}
   +b\cos(\theta)\frac{\partial\mu}{\partial x} & = 0.
  \endalign
   $$
By the Frobenius Theorem, or equivalently checking that mixed derivatives
are equal, we see that this system has a solution if and only if
$\,b\sin(\theta)\frac{\partial\mu}{\partial x}=0$. 
The second equation then implies that $\sigma$ only depends on $\theta$,
so the first equation reduces to an ordinary differential equation.
The general solution to this system of equations is, therefore,
$$ \sigma=\sigma(\theta),\quad
   \mu=\mu_0\cos(\theta)-\frac{1}{b}\cos(\theta)
   \int\frac{\partial\sigma}{\partial\theta}\sec^2\theta\,d\theta.
   $$
For the rest of this example, we will use the particular solution:
$$ \sigma=\sigma_0,\quad \mu=\mu_0\cos(\theta).
   \tag15 $$

It is easiest to find $\,\widehat g\,$ using the basis
$\frac{\partial}{\partial\theta}$, 
$\,\lambda\left(\frac{\partial}{\partial
\theta}\right)$. 
By the relation $\,\widehat g\lambda=g$,
we already know
$\,\widehat g\left(\lambda(\frac{\partial}{\partial\theta}),
\frac{\partial}{\partial\theta}\right)\,$ 
and
$\,\widehat g\left( \lambda\left(\frac{\partial}{\partial\theta}\right)
 \lambda\left(\frac{\partial}{\partial\theta}\right)\right)$,
so we only need to find $\,\widehat g\left(
\frac{\partial}{\partial\theta},\frac{\partial}{\partial\theta}\right)$.
The Lie bracket,
$$
\left[ \frac{\partial}{\partial\theta},
 \lambda(\frac{\partial}{\partial\theta})\right]=-\mu_0\sin(\theta)
 \frac{\partial}{\partial x}=\sigma_0\tan(\theta)
 \frac{\partial}{\partial\theta}-\tan(\theta)\lambda
 \left( \frac{\partial}{\partial\theta}\right) 
 $$
may be plugged into the $\widehat g$-Equation (10) to get:
$$ 
\sigma_0\frac{\partial \widehat g_{11}}{\partial\theta}
   +\mu_0\cos(\theta)\frac{\partial \widehat g_{11}}{\partial x}
   +2\sigma_0\tan(\theta)\widehat g_{11}
   -2\tan(\theta)=0.
   $$
The surface $\theta=0$ is non-characteristic, so the above
equation may be solved with initial data 
$\;\widehat g_{11}(0,x)=h(x)$.

The flow of the vector field $\lambda PX$ is given by:
$$ 
\dot\theta=\sigma_0, \quad \dot x=\mu_0\cos(\theta),
  $$
or
$$ 
\theta=\sigma_0t,\quad x=a+\frac{\mu_0}{\sigma_0}\sin(\sigma_0t), 
  $$
using $\theta=0$ at $t=0$. With this flow, the $\widehat g$ equation
becomes:
$$ 
\frac{d\widehat g_{11}}{dt} +2\sigma_0\tan(\sigma_0t)
   \widehat g_{11}-2\tan(\sigma_0t)=0,
   $$
so
$$ 
   \widehat g_{11} = \frac{1}{\sigma_0}
     +\left(h(a)-\frac{1}{\sigma_0}\right) \cos^2(\sigma_0t)
   =\frac{1}{\sigma_0}+
       \left(h\left(x-\frac{\mu_0}{\sigma_0}
       \sin(\theta)\right)-\frac{1}{\sigma_0}\right)\cos^2(\theta).
$$
Again, we only take a particular solution:
$$ 
\widehat g_{11}=\frac{1}{\sigma_0}+r\cos^2(\theta),  
$$
where $r$ is a constant.
The equations
$$ 
 b\cos(\theta)=
   g\left( \frac{\partial}{\partial\theta},\frac{\partial}{\partial x}\right)
  =\widehat g  \left( \lambda
    (\frac{\partial}{\partial\theta}),\frac{\partial}{\partial x}\right)
  =\sigma_0\widehat g_{12}+\mu_0\cos(\theta)\widehat g_{22}
$$
and 
$$
1=g\left( \frac{\partial}{\partial\theta},
     \frac{\partial}{\partial\theta}\right)
  =\widehat g  \left( \lambda
    (\frac{\partial}{\partial\theta}),
   \frac{\partial}{\partial\theta}\right)
  =\sigma_0\widehat g_{11}+\mu_0\cos(\theta)\widehat g_{12} 
$$
may now be solved for the remaining terms of $\,\widehat g\,$:
$$ 
\widehat g_{12}=-\frac{\sigma_0}{\mu_0}r\cos(\theta),
   \qquad \widehat g_{22}=\frac{b}{\mu_0}+\frac{\sigma^2_0}{\mu^2_0}r.
   $$

We can compute $\widehat V$ in the same way. Using the $\lambda XP$ flow
and initial data $\widehat V(0,x)=w(x)$,
the $\,\widehat V$-Equation (11) may be written as an ordinary differential
equation:
$$ 
\frac{d\widehat V}{dt}=-\sin(\sigma_0t). 
$$
Thus,
$$    
\widehat V  =\frac{1}{\sigma_0}(\cos(\sigma_0t)-1)+w(a)\,
  =\,\frac{1}{\sigma_0}(\cos(\theta)-1)
+w\left( x-\frac{\mu_0}{\sigma_0}\sin\theta\right).
$$
Finally, solving for $\,\widehat c$, we get:
$$ 
\widehat c(X)=K(\theta,x,\dot\theta,\dot x)
   \left( b\cos(\theta)\frac{\partial}{\partial\theta}
     -\frac{\partial}{\partial x}\right)
   $$
with an arbitrary function $\,K$.

The functions $\widehat g$, $\widehat V$ and $\widehat c$ produce a 
family of control
laws via equation (3). In this example, the control laws depend on
four functions, $\,\sigma_0(\theta)$, $h(x)$, $w(x)$, and 
$\,K(\theta,x,\dot\theta,\dot x)$, and one
constant, $\mu_0$. We have already chosen $\sigma_0$ and $h$ to
be constants. We will address the question of how to choose the
unknown functions and parameters in order to best meet specific
design criteria in a future paper. For now, we will just pick elementary
functions that will insure that $\theta=0$, $x=0$ is an asymptotically
stable equilibrium. If $(0,0)$ is an equilibrium, 
$\,\widehat  {grad}\,\widehat V (0,0)$
must be zero, so $w'(0)$ must be zero.
One standard way to insure that $(0,0)$ will be an
asymptotically stable equilibrium is to pick positive definite
$\widehat g$, $\widehat V$ and $\widehat g\widehat c$, then
$ \widehat H (X)\equiv\frac12 \widehat g(X,X)+\widehat V$ will be a
Lyapunov function with time rate of change 
$\,-\widehat g\,(\widehat c(X),X)$.
The Hessian of $\widehat V$ at $(0,0)$ is
$$ 
D^2_{(0,0)} =
  \left(\aligned \frac{\mu^2_0}{\sigma^2_0}w''(0)-\frac{1}{\sigma_0}
    & \quad -\frac{\mu_0}{\sigma_0}w''(0) \\
    -\frac{\mu_0}{\sigma_0}w''(0)\quad\qquad &\qquad w''(0)
    \endaligned \right).
  $$
So we should require:
$$ 
w''(0)>0\quad\text{and}\quad
   \det D^2_{(0,0)} \widehat V=-\frac{1}{\sigma_0}w''(0)>0.
   $$
We will choose $\sigma_0<0$, and $\,w(z)=\frac12 w_1z^2\,$ with $\,w_1>0$.

The mass matrix should also be positive definite.
Thus we should have $\widehat g_{11}>0$ and $\widehat g_{22}>0$ so that $r>0$
and $b\mu_0+\sigma^2_0r>0$. In addition,
the determinant of the mass matrix will also be positive:
$$ 
\frac{b}{\sigma_0\mu_0} + \frac{br}{\mu_0}\cos^2(\theta)
    +\frac{\sigma_0}{\mu^2_0}r>0.
    $$
Rearranging and  using the fact that $\sigma_0<0$ gives:
$$ -\sigma_0b\mu_0 r\cos^2(\theta)
   >\sigma^2_0r+b\mu_0.
   $$
So,
$$ 
\mu_0>\frac{\sigma^2_0r+b\mu_0}{-\sigma_0br\cos^2(\theta)}>0,
   $$
and
$$ 
\cos^2(\theta)>\frac{\sigma^2_0r+b\mu_0}{-\sigma_0b\mu_0r}.
    \tag16 
$$
Assuming $\mu_0>0$, $r>0$ and $\sigma_0<0$,
the mass matrix will be positive definite provided condition (16)
holds. When the mass matrix and potential energy of the model system
are positive definite, the point $(0,0)$ will be Lyapunov stable for any
positive
semi-definite $\widehat g\widehat c$. In particular,
when $\,K\equiv 0\,$ (i.e., $\,\widehat c\equiv 0$), 
we will get Lyapunov stability. 
We call controllers
with $\,\widehat c\equiv 0\,$ 
conservative controllers.
A direct computation shows that:
$$ 
\align   \widehat g(\widehat c(X),X)
  & =K\widehat g \left( b\cos(\theta)
     \frac{\partial}{\partial\theta}-\frac{\partial}{\partial x},
     \dot\theta\frac{\partial}{\partial\theta}
     +\dot x\frac{\partial}{\partial x} \right)     \\
  &=K \left[(b\cos(\theta)\widehat g_{11}-\widehat g_{12})\dot\theta
     +(b\cos(\theta)\widehat g_{12}-\widehat g_{22})\dot x\right] \\
  &=(\det \widehat g)\cdot K\cdot (\mu_0\cos\theta\dot\theta-\sigma_0\dot x).
  \endalign 
$$                     
This will be positive semi-definite if
$K=\Phi(\theta,x)(\mu_0\cos\theta\dot\theta-\sigma_0\dot x)$, 
where $\,\Phi\,$ is any positive function. (We can do no better because
$\lambda PX=\sigma_0\frac{\partial}{\partial\theta}+\mu_0\cos\theta
 \frac{\partial}{\partial x}$
will be a zero mode of any
admissible $\,\widehat c$.)  
We will see that this is sufficient to prove that $(0,0)$ is an
asymptotically stable equilibrium of the controlled system.
Equation (3), equation (7) and our values for $\widehat g$, 
$\widehat V$ and $\widehat c$
combine together to give:
$$
\aligned
 & u=g\left( f,\frac{\partial}{\partial x}\right)
    =\left( b+\frac{r\det  g }{\mu_0\det \widehat g }\right)
     \left(\cos(\theta)\sin(\theta)-\sin(\theta)\dot\theta^2 \right) \\
 &\qquad  -\frac{w_1\det  g }{\sigma_0\det \widehat g }
     \left(x-\frac{\mu_0}{\sigma_0}\sin(\theta)\right)
    +\det g\; \Phi(\theta,x)\,(\mu_0\cos(\theta)\dot\theta-\sigma_0\dot x),
  \endaligned 
\tag 17
$$
where
$$
\det g=1-b^2\cos^2(\theta),\qquad
\det \widehat g=\frac{b}{\sigma_0\mu_0}
+\frac{br}{\mu_0}\cos^2(\theta)
 +\frac{\sigma_0 r}{\mu_0^2}. \tag 18
$$
Thus, we have almost proved the following result.


\proclaim{Proposition 5} Let constants $\,\mu_0$, $\,\sigma_0$, $\,r$, 
and $\,w_1\,$ satisfy the conditions
$$
\mu_0>0,\quad\sigma_0<0,\quad w_1>0,\quad 
1\,>\,{\sigma^2_0 r+b\mu_0\over -\sigma_0 b\mu_0 r}. \tag 19
$$
Let $\,\Phi(\theta,x)\,$ be any strictly positive function. 

Then $\,(0,0)\,$ is an asymptotically stable equilibrium 
of the controlled system (14) with the control law defined 
by (17), (18).
\endproclaim
\demo{Proof} Define the controlled Hamiltonian 
$\widehat H =\frac12\widehat g(\dot\gamma,\dot\gamma)+\widehat V$. 
We have shown above that the Hessian of $\,\widehat H\,$ 
is positive definite in some neighborhood of $\,(0,0)\,$ 
when conditions (19) hold. 
The time derivative of $\widehat H $ is:
$$ 
\align  \frac{d}{dt}(\widehat H )
  &=\frac{d}{dt} \left( \frac12 \widehat g (\dot\gamma,\dot\gamma)
     +\widehat V \right) \\
  &=\widehat g (\nabla_{\dot\gamma}\dot\gamma 
       +\ \widehat  {grad}\widehat V,\dot\gamma)\\
  &=-\widehat g (\widehat c(\dot\gamma),\dot\gamma) \\
  &=-\det \widehat g \;\Phi(\theta,x)\,
(\mu_0\cos(\theta)\dot\theta-\sigma_0\dot x)^2.
  \endalign 
$$ 
Check that there is no solution to the controlled equations satisfying
$\mu_0\cos(\theta)\dot\theta - \sigma_0 \dot x\equiv 0$. 
Since $\,\Phi>0\,$ by assumption 
and $\,\det\widehat g\big|_{(\theta,x)}\,>0\,$ for all 
$\,(\theta,x)\,$ sufficiently close to $\,(0,0)$, 
the Hamiltonian $\,\widehat H\,$ decreases along the solutions of 
(17), and, therefore, may serve as a Lypunov function.  
The (local) 
asymptotic stability follows from Lyapunov's Theorem.
\qed
\enddemo
 
\medskip


\proclaim{Remark 1} 
It is not possible to construct a control law so that $(0,0)$ is a
globally asymptotically stable equilibrium. 
\endproclaim

Indeed, the solutions of the
controlled system for such a control law would produce a continuous
function,
$F:[0,\infty]\times TQ\to TQ$ so that $F_0$ is the identity
map and $F_\infty$ is the constant map. In other words, the
flow of a vector field with a globally asymptotically stable
equilibrium is a contraction. However, for the inverted pendulum cart,
$TQ\simeq S^1\times{\Bbb R}^3$, which is not contractible.
\medskip


Since the inverted pendulum cart cannot be globally stabilized
by a control law, we can only try to maximize the size of the
basin of attraction. We will
compare a special case of our nonlinear control law with a linear
control law that has been implemented on an inverted pendulum
cart in our lab. Our cart has
$$ 
\align
   M=5.02\ \text{Kg} & \qquad m=.454\ \text{Kg} \\
   \ell=.425\ m\  & \qquad\ I=.11\ \text{Kg}\ m^2
   \endalign 
$$
giving $b=.188$. The linear control law that we obtained by placing poles
at $-5$, $-6$, and a double pole at $-2$ is:
$$
u_{\text{lin}} = 1021 \theta + 115.8 x + 918.5 \dot\theta +158.2 \dot x \,.
$$
This is only an approximation of the actual control law that is running
in the lab. In the lab we had to model the DC motor, the observer, and use
a discrete control law.

Thus far, we have chosen each of the arbitrary functions in our non-linear
control law to be a constant, and we took 
$\,w(x)=\frac{1}{2} w_1 x^2$, 
since it has to have positive second derivative.
As a rough guess we took $\mu_0 = 10$, $\sigma_0 = -.05$, $r=1000$,
and $w_1 = 1.5$. We also set the function $\Phi = 1$.  These constants 
satisfy  conditions (19). We chose the functions
to be constants when possible just to simplify the exposition.  

\proclaim{ Remark 2 } This would
not be the best choice for engineering applications. In particular, 
$-\sigma_0$ controls the coefficient of $\dot x$ in our control law. 
In order
to stabilize the system after a large angular disturbance, we would like
$-\sigma_0$ to be small so that the cart will be free to accelerate
under the pendulum. On the other hand, $-\sigma_0$ needs to be much bigger
in order to insure that the time constant is reasonable. 
We hope to describe some practical approches to choosing the arbitrary
functions which appear in our family of control laws so that the resulting
controler will meet specific design criteria in a future paper.
\endproclaim
\medskip

For a first rough test of our control law, we ran numerical simulations
of the system using both our control law and a linear control law with
the arbitrary functions and parameters specified above. We used 10000
diferent sets of initial conditions.  It appears that the linear controler
has a shorter settling time than the nonlinear control law, when it does
not blow up. However, the nonlinear control law appears to stabilize the
system starting from any of the initial conditions stabilized by the
linear control law, and also appears to stabilize the system starting from
many initial conditions for which the linear control law 
produces an unbounded response.
\bigskip
\settabs 2 \columns
\+ \hfill\psfig{file=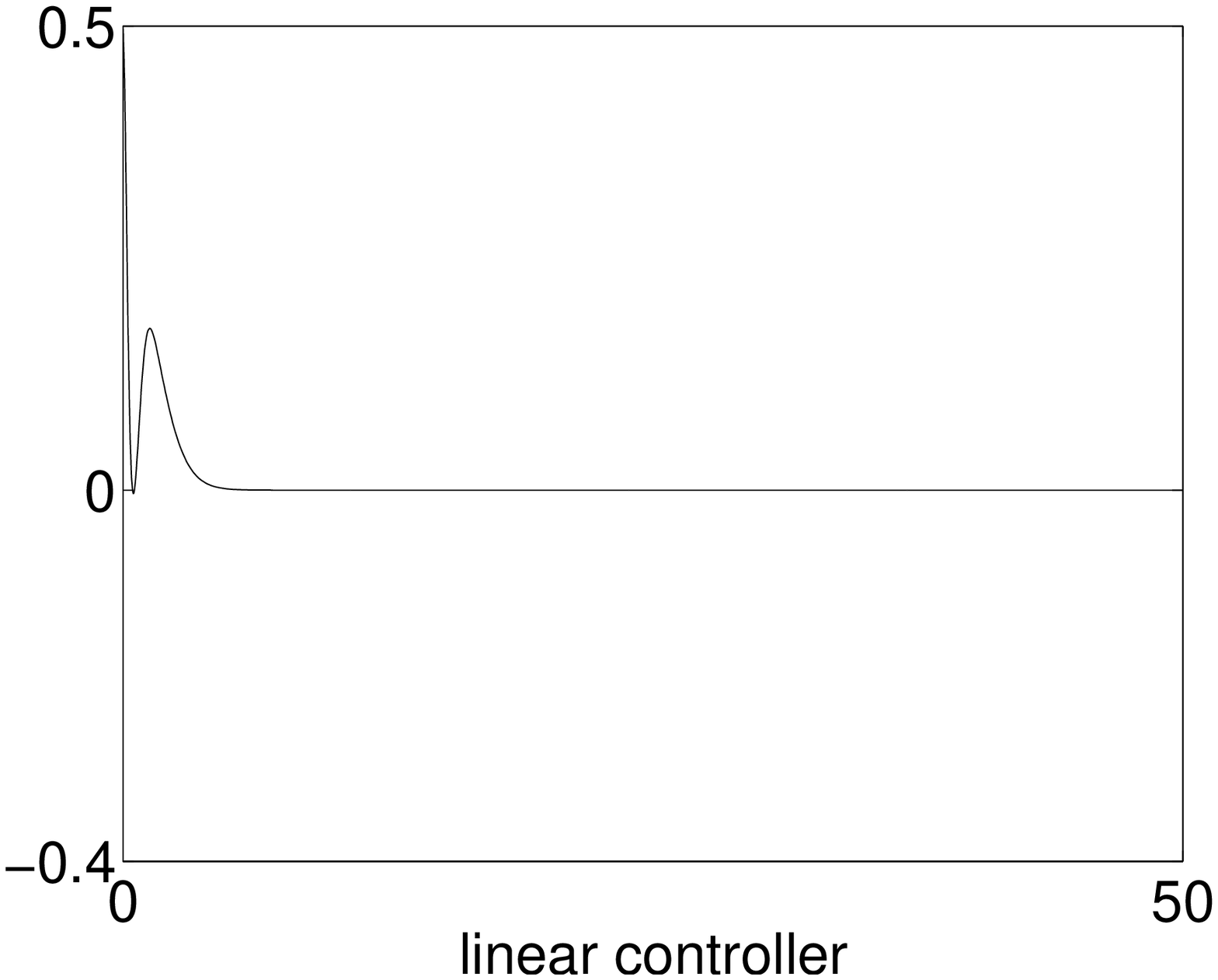,width=2.2truein,height=2.2truein}\hfill & 
   \hskip40bp\psfig{file=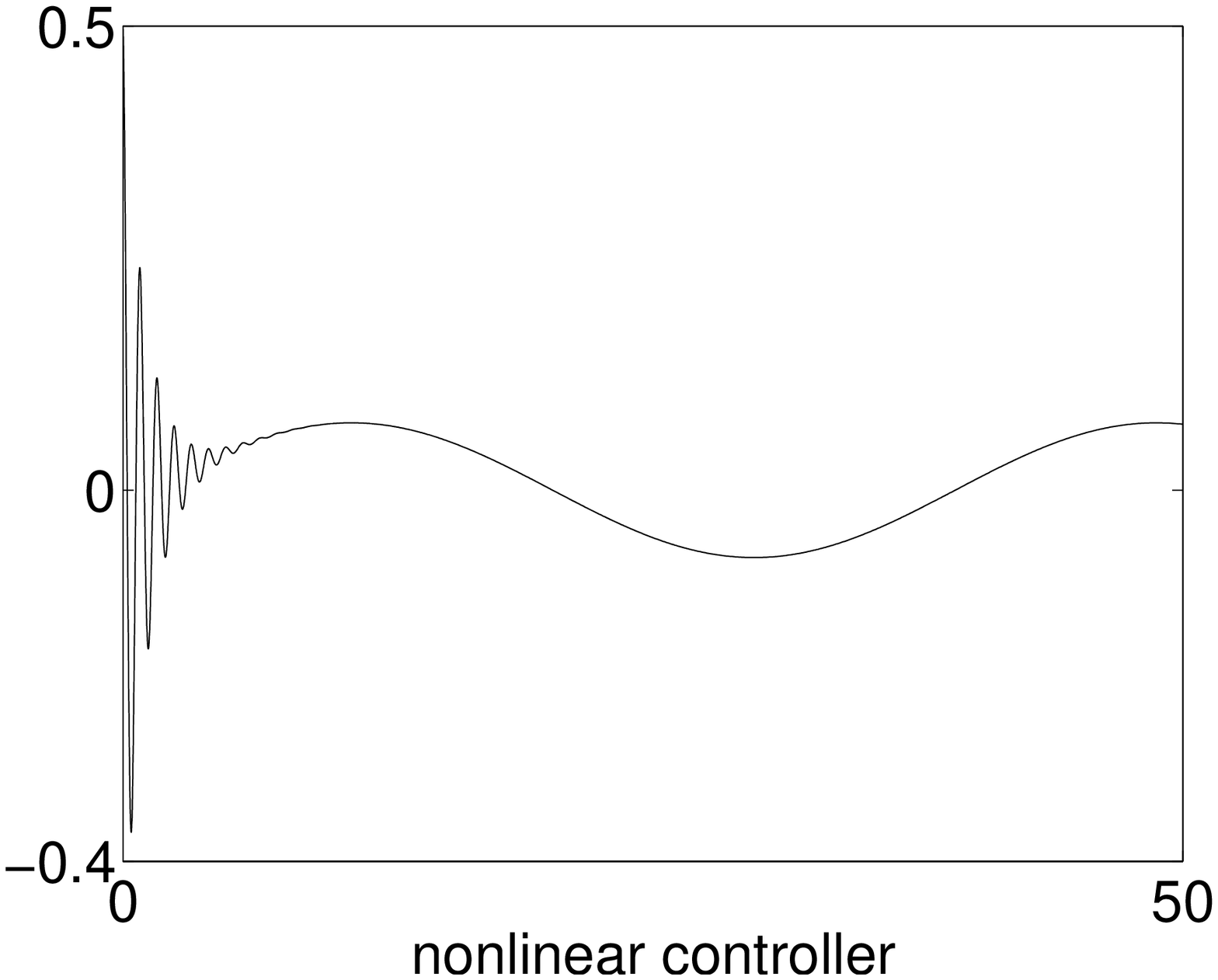,width=2.2truein,height=2.2truein}\hfill\cr
\centerline{Figure 2}

Figure 2 shows graphs of the angular position versus time produced using
the linear control law and the nonlinear control law starting from initial
conditions: $\theta = .5$, $\dot\theta =-.5$, $x=0$, and
$\dot x=0$. The graphs of the other state variables and control input are
qualitatively similar. This output is typical of the responce obtained when
both control laws appear to stabilize the system. 
\bigskip
\settabs 2 \columns
\+ \hfill\psfig{file=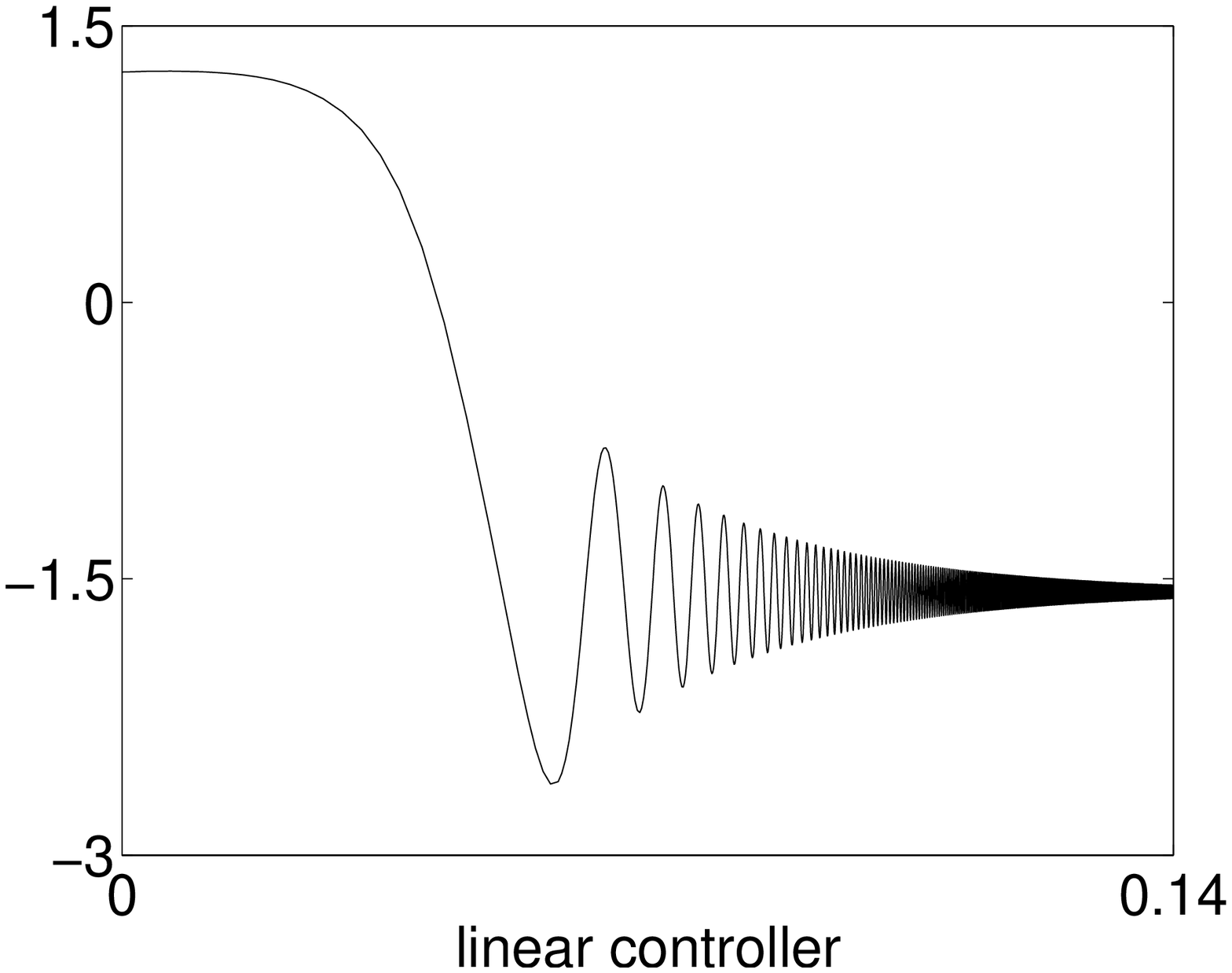,width=2.2truein,height=2.2truein}\hfill & 
   \hskip40bp\psfig{file=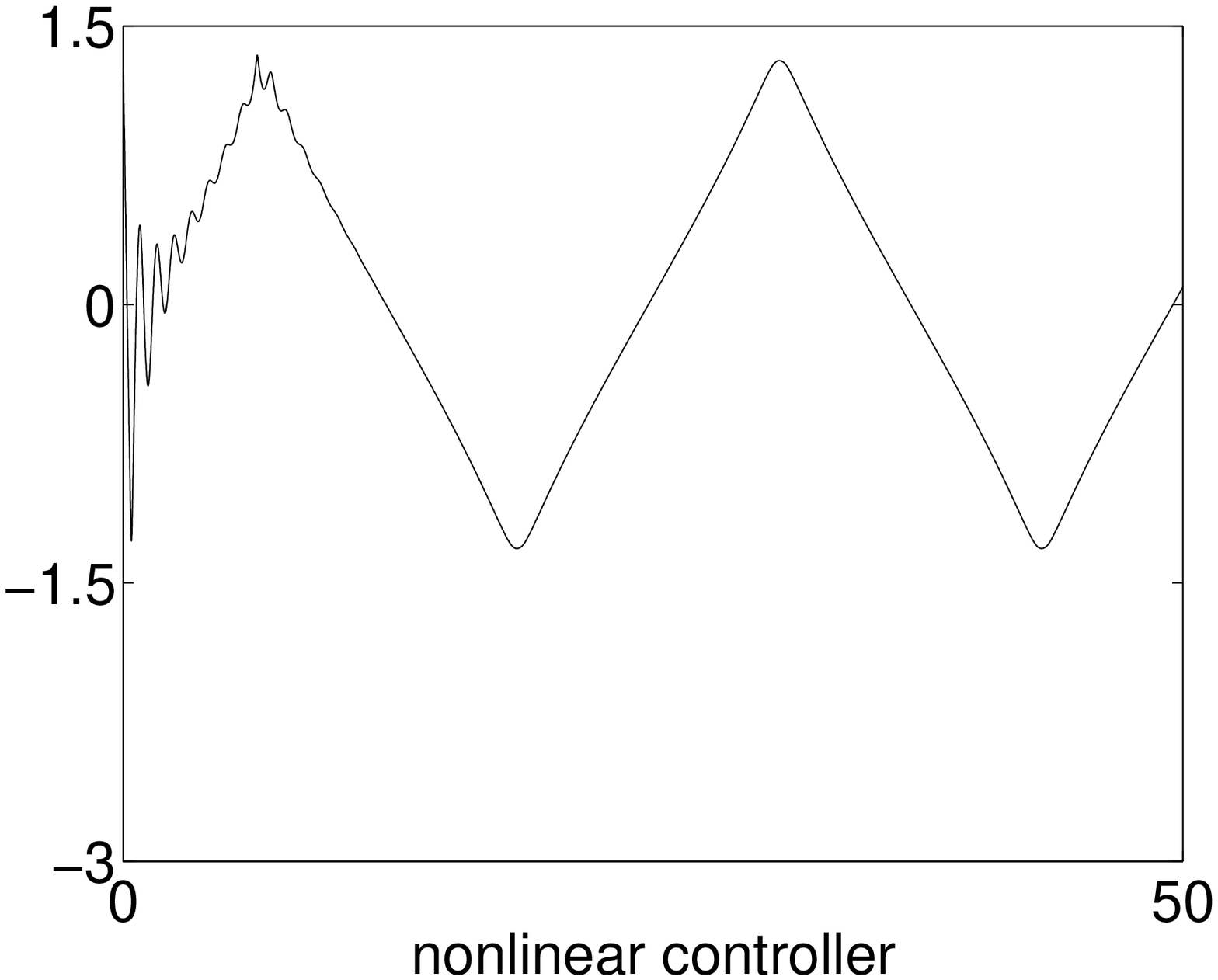,width=2.2truein,height=2.2truein}\hfill\cr
\centerline{Figure 3}

Figure 3 shows graphs of the angular position versus time produced using
the linear control law and the nonlinear control law starting from initial
conditions: $\theta = 1.25$, $\dot\theta =1.3$, $x=0$, and
$\dot x=0$. For the same initial conditions, Figure 4 
shows graphs of the cart position versus time produced using
the linear control law  
and the nonlinear control law. 

\medskip
\settabs 2 \columns
\+ \hfill\psfig{file=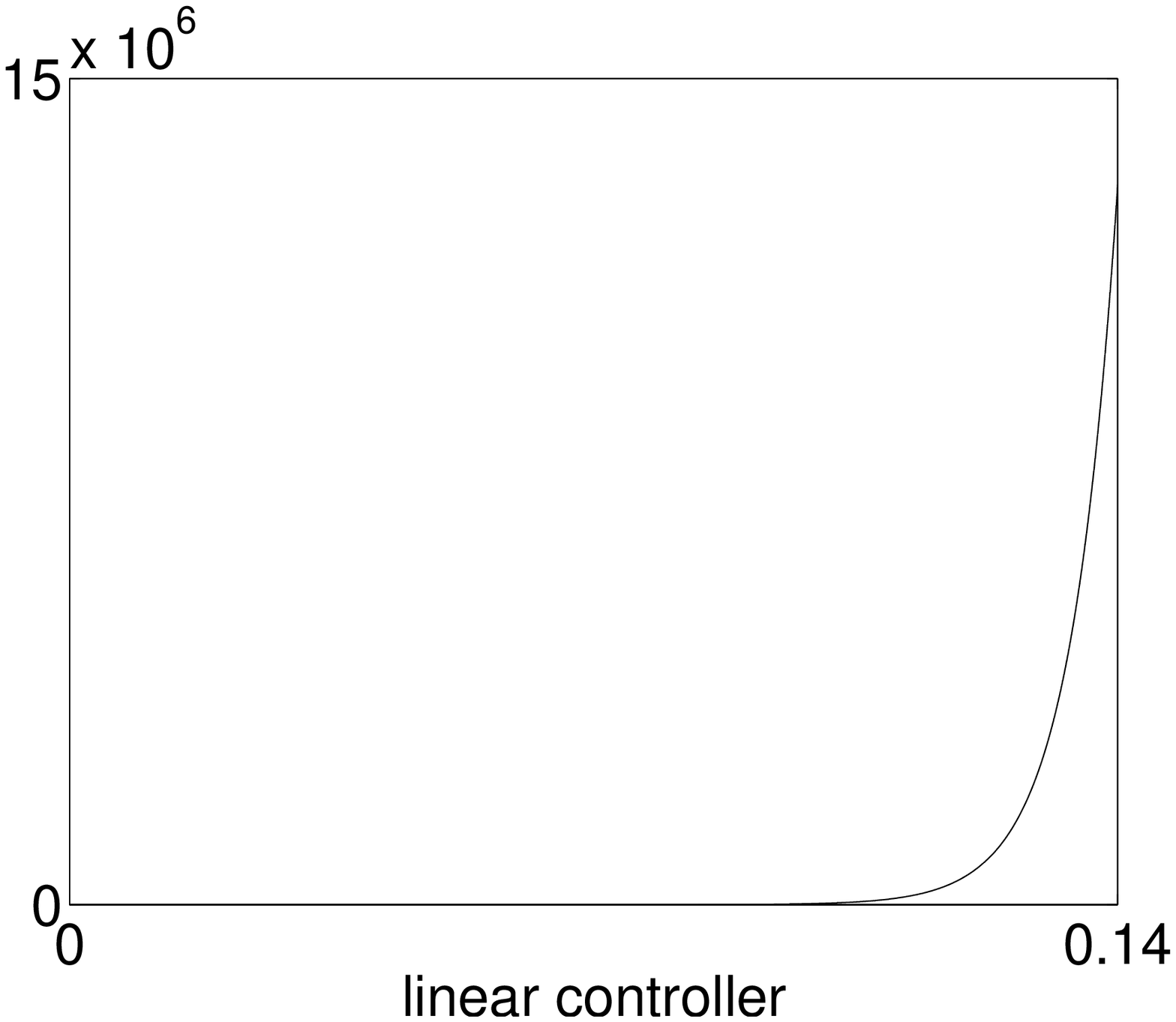,width=2.2truein,height=2.2truein}\hfill & 
   \hskip40bp\psfig{file=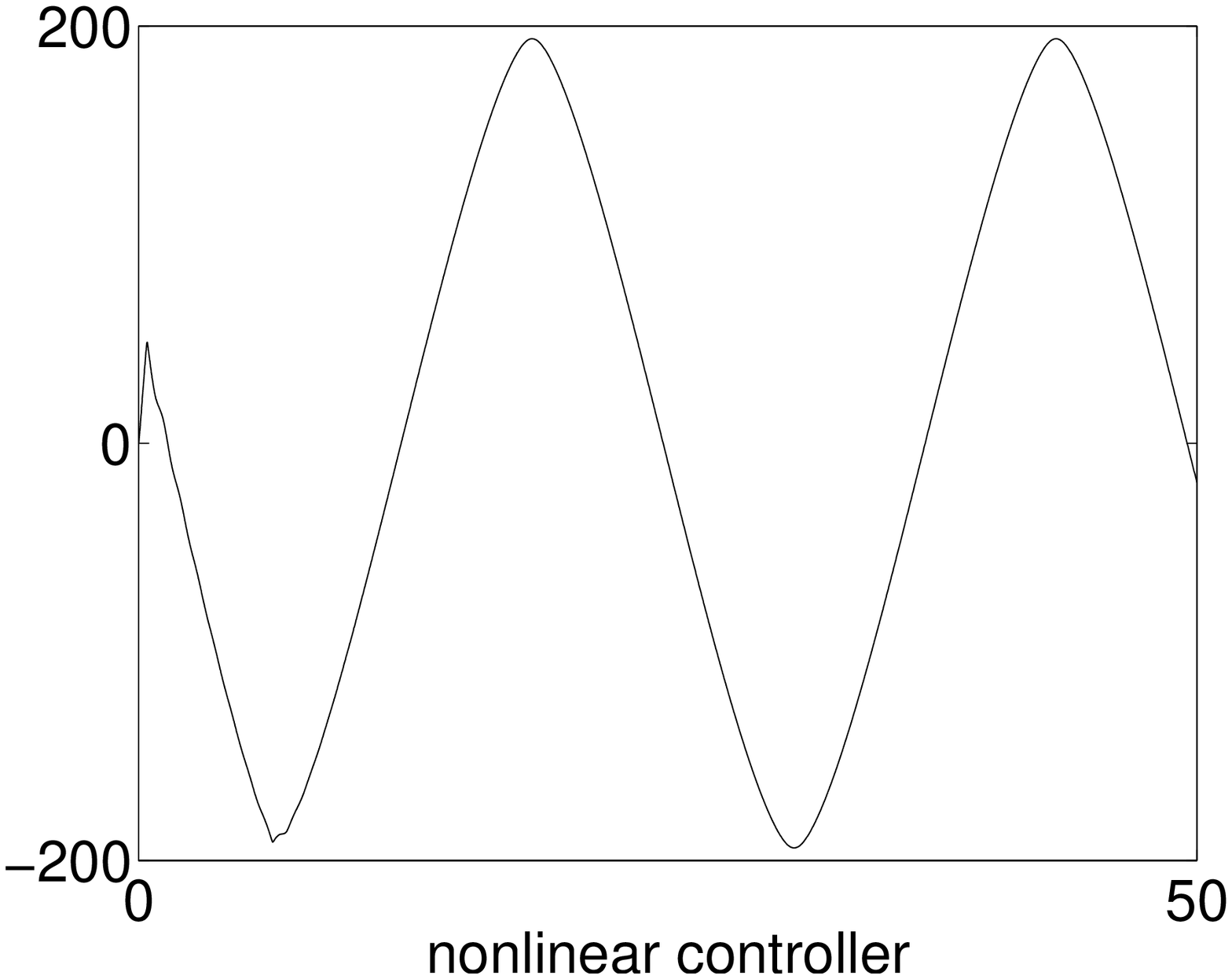,width=2.2truein,height=2.1truein}\hfill\cr
\centerline{Figure 4}

\noindent For the linear controller, $\,x(t)\,$ increases exponentially. 
For the nonlinear controller, $\,x(t)\,$ stabilizes to $0$, but 
the system is underdamped, as we discussed previously.
The graphs of the other state variables and control input are
again qualitatively similar. These graphs are representative for initial
conditions for which the linear control blows up, and the nonlinear control
stabilizes.  All of our numerical simulations were produced using MATLAB 5.1.
\medskip

Our infinite dimensional family of control laws contains a one
dimensional family of control laws, which was found previously.
Bloch, Leonard and Marsden have developed a method for constructing
control laws for mechanical systems with symmetry, 
\cite{Bloch, Lenard,  
Marsden, 1997}. Their
control law will retain the symmetry, so it cannot produce a truly
Lyapunov stable equilibrium. It will only produce a Lyapunov
stable equilibrium in shape space. This means that their method
will not work for an inverted pendulum cart on a hill, say. If the
cart is on level ground, the quantity $\dot x+A\cos(\theta)\dot\theta$
will be
conserved. The quantity $\dot x+A\cos(\theta)\dot\theta$
is conserved in the closed-loop dynamics of any of the control laws
developed by Bloch, Leonard and Marsden. It follows that
any solution which passes through a point with $\dot\theta=0$
and $\dot x\not= 0$ will run off to infinity.
\medskip

\head LINEARLY NON-CONTROLLABLE EXAMPLE \endhead

The final example that we will consider is an abstract system which is
open-loop unstable with the property that the linearized system
cannot be stabilized by any control law.
Our method will generate a control law with a globally asymptotically
stable equilibrium. 

The abstract system is:
$g=dx^2+dy^2$, $P=dy\otimes \frac{\partial}{\partial y}$ and
$V=-\frac32 x^4+45x^2y^2+32xy^3$.
In coordinates, the control problem reads:
$$ 
\align
  & \ddot x-6x^3+90xy^2+32y^3=0 \\
  &\ddot y+90x^2y+96xy^2=u.
  \endalign 
$$
This system is, clearly, open-loop $(u=0)$ unstable. For example,
$x=(\varepsilon^{-1}-\sqrt{3}t)^{-1}$, $y=0$ is a solution to the
system for any positive $\varepsilon$. The linearized version of
the system is:
$\;\ddot x=0$, $\ddot y=u$. Clearly, $u$ can have no effect on $x$.

To apply our method to this system, let
$\lambda PX=\sigma\frac{\partial}{\partial x}
  +\mu\frac{\partial}{\partial y}$,
then the $\lambda$-Equation (9) reads:
$$  
\frac{\partial\sigma}{\partial x}=0,\qquad\qquad
    \frac{\partial\sigma}{\partial y}=0.
    $$
Thus, $\,\sigma$ is any constant and $\mu$ is any function.
Pick $\sigma=1$ and $\mu=1$. The $\,\widehat g$-Equation (10) becomes:
$$ 
\frac{\partial \widehat g_{11}}{\partial x} 
+ \frac{\partial \widehat g_{11}}{\partial y}
   =0. 
$$
Pick $\widehat g_{11}=2$. The equations $\widehat g\lambda=g$ lead to
$\widehat g_{22}=1$ and $\widehat g_{12}=-1$.
It is easy to check that the model mass matrix is positive definite.
The flow equation for $\widehat V$ is
$$ \frac{\partial \widehat V}{\partial x}
   + \frac{\partial \widehat V}{\partial y}
   =-6x^3+90xy^2+32y^3.
   $$
It is not hard to check that
$$ 
\widehat V=(x^2-3xy)^2+(x^2-4xy-2y^2)^2 
$$
is a positive definite solution to this equation. Finally,
pick $\widehat c=(\dot y-\dot x)\frac{\partial}{\partial y}$ as a solution to
$P(\widehat c)=P(c)=0$. As before, we show that there is no
non-trivial solution to the controlled equations with
$\dot y-\dot x\equiv 0$. Thus, $(0,0)$ is a 
globally asymptotically stable equilibrium.

\bigskip


\centerline {REFERENCES}
\medskip

\ref\by A. Bloch, N. Lenard, J. Marsden
\paper Stabilization of Mechanical Systems Using Controlled Lagrangians
\jour Proc. 1997 IEEE Conference on Decision and Control
\yr 1997 \pages 2356--2361
\endref

\ref \by N. Hicks
\book Notes on Differential Geometry
\publ Van Nostrad \yr 1965
\endref


\end


Figure 3 shows the graph of the angular responce starting
at the above initial conditions over a longer time scale using the nonlinear
control law.

\psfig{file=lc125_13.eps,width=2in,height=2in} & 
\psfig{file=t_nl_ex.eps,width=2in,height=2in}\cr

\psfig{file=lc125_13.eps,width=2.5truein,height=2.5truein} &
\psfig{file=t_nl_ex.eps,width=2.5truein,height=2.5truein}\cr

\vfill
\pagebreak
\quad
\vskip 6.5in
\centerline{Figure 2}

\vfill
\pagebreak
\quad
\vskip 6.5in
\centerline{Figure 3}
\vfill
\pagebreak
\quad
\vskip 6.5in
\centerline{Figure 4}

\vfill
\pagebreak

\vskip.1in

\vskip.1in

, \cite{John, 1982}

\ref\by {\bf [HS]\/}\ \ \ M. Hirsch, S. Smale
\book Differential Equations, Dynamical Systems,
and Linear Algebra
\publ Academic Press \yr 1974
\endref

\ref\by {\bf [J]\/}\ \ \ \ \ \ \ F. John
\book Partial Differential Equations
\publ Springer Verlag \yr 1982
\endref